\documentclass{article}%
\usepackage{amsmath}
\usepackage{amsfonts}
\usepackage{amssymb}
\usepackage{graphicx}%
\providecommand{\U}[1]{\protect\rule{.1in}{.1in}}
\begin{document}

\title{Inferences from Prior-based Loss Functions}
\author{Michael\ Evans\\Department of Statistics\\University of Toronto\\Toronto, ON M5S 3G3\\mevans@utstat.utoronto.ca
\and Gun Ho Jang\\Department of Biostatistics and Epidemiology\\University of Pennsylvania\\Philadelphia, PA 19104, USA\\gunjang@upenn.edu}
\date{}
\maketitle

\begin{abstract}
Inferences that arise from loss functions determined by the prior are
considered and it is shown that these lead to limiting Bayes rules that are
closely connected with likelihood. The procedures obtained via these loss
functions are invariant under reparameterizations and are Bayesian unbiased or
limits of Bayesian unbiased inferences. These inferences serve as
well-supported alternatives to MAP-based inferences.

\end{abstract}


\section{Introduction}

Suppose we have a sampling model, given by a collection of densities
$\{f_{\theta}:\theta\in\Theta\}$ with respect to a support measure $\mu$ on
sample space $\mathcal{X},$ and a proper prior, given by density $\pi$ with
respect to support measure $\nu$ on $\Theta.$ When we observe data $x$ these
ingredients lead to the posterior on $\Theta$ with density given by
$\pi(\theta\,|\,x)=\pi(\theta)f_{\theta}(x)/m(x)$ with respect to support
measure $\nu$ where $m(x)=\int_{\Theta}\pi(\theta)f_{\theta}(x)\,\nu
(d\theta).$

One can determine inferences based on these ingredients alone. For example,
suppose we are interested in a characteristic $\psi=\Psi(\theta)$ where
$\Psi:\Theta\rightarrow\Psi$ and we let $\Psi$ stand for both the space and
mapping to conserve notation. The highest posterior density (hpd), or
MAP-based, approach to determining inferences constructs credible regions of
the form
\begin{equation}
H_{\gamma}(x)=\{\psi_{0}:\pi_{\Psi}(\psi_{0}\,|\,x)\geq h_{\gamma}(x)\}
\label{eq1}%
\end{equation}
where $\pi_{\Psi}(\cdot\,|\,x)$ is the marginal posterior density with respect
to a support measure $\nu_{\Psi}$ on $\Psi,$ and $h_{\gamma}(x)$ is chosen so
that $h_{\gamma}(x)=\sup\{k:\Pi_{\Psi}(\{\psi:\pi_{\Psi}(\psi\,|\,x)\geq
k\}\,|\,x)\geq\gamma\}.$ It follows from (\ref{eq1}) that, if we want to
assess the hypothesis $H_{0}:\Psi(\theta)=\psi_{0},$ then we can use the tail
probability given by $1-\inf\{\gamma:\psi_{0}\in H_{\gamma}(x)\}.$
Furthermore, the class of sets $H_{\gamma}(x)$ is naturally "centered" at the
posterior mode (when it exists uniquely) as $H_{\gamma}(x)$ converges to this
point as $\gamma\rightarrow0.$ The use of the posterior mode as an estimator
is commonly referred to as MAP (maximum \textit{a posteriori}) estimation. We
can then think of the size of the set $H_{\gamma}(x),$ say for $\gamma=0.95,$
as a measure of how accurate the MAP estimator is in a given context.
Furthermore, we have that when $\Theta$ is an open subset of a Euclidean
space, then $H_{\gamma}(x)$ minimizes volume among all $\gamma$-credible
regions. The use of MAP-based inferences is very common in machine learning
contexts, see, for example, Bishop (2006).

It is well-known, however, that hpd inferences suffer from a serious defect.
In particular, in the continuous case hpd inferences are not invariant under
reparameterizations. For example, this means that if $\psi_{\text{MAP}}(x)$ is
the MAP estimate of $\psi$, then it is not necessarily true that
$\Upsilon(\psi_{\text{MAP}}(x))$ is the MAP estimate of $\tau=\Upsilon(\psi)$
when $\Upsilon$ is a 1-1, smooth transformation. The noninvariance of a
statistical procedure seems very unnatural as it implies that the statistical
analysis depends on the parameterization and typically there does not seem to
be a good reason for this.

A class of inferences, similar to hpd inferences, avoids this lack of
invariance. These are referred to as \textit{relative surprise inferences} and
are based on the regions
\begin{equation}
C_{\gamma}(x)=\{\psi:\pi_{\Psi}(\psi\,|\,x)/\pi_{\Psi}(\psi)\geq c_{\gamma
}(x)\} \label{eq2}%
\end{equation}
where $\pi_{\Psi}$ is the marginal prior density with respect to a support
measure $\nu_{\Psi}$ on $\Psi,$ and $c_{\gamma}(x)=\sup\{k:\Pi_{\Psi}%
(\{\psi:\pi_{\Psi}(\psi\,|\,x)/\pi_{\Psi}(\psi)\geq k\}\,|\,x)\geq\gamma\}.$
The hypothesis $H_{0}:\Psi(\theta)=\psi_{0}$ is assessed by computing the tail
probability
\begin{equation}
1-\inf\{\gamma:\psi_{0}\in C_{\gamma}(x)\}=\Pi_{\Psi}(\pi_{\Psi}%
(\psi\,|\,x)/\pi_{\Psi}(\psi)\leq\pi_{\Psi}(\psi_{0}\,|\,x)/\pi_{\Psi}%
(\psi_{0})\,|\,x). \label{rspvalue}%
\end{equation}
We refer to $\pi_{\Psi}(\psi\,|\,x)/\pi_{\Psi}(\psi)$ as the \textit{relative
belief ratio} of $\psi$ as it measures how beliefs in $\psi$ being the true
value change from \textit{a priori} to \textit{a posteriori}. The relative
surprise terminology then comes from (\ref{rspvalue}) as this is measuring how
surprising the value $\psi_{0}$ is by comparing its relative belief ratio to
the relative belief ratios of other values of $\psi.$ The corresponding
estimator is given by the maximizer of the ratio $\pi_{\Psi}(\psi
\,|\,x)/\pi_{\Psi}(\psi),$ which we refer to as the \textit{least relative
surprise estimator} (LRSE), and denote as $\psi_{\text{LRSE}}(x).$ Note that
$\psi_{\text{LRSE}}(x)$ is the least surprising value as it maximizes
(\ref{rspvalue}). Beyond their invariance these inferences have many
optimality properties in the class of all Bayesian inferences as documented in
Evans (1997), Evans, Guttman and Swartz (2006), Evans and Shakhatreh (2008)
and Jang (2010). In this paper we will establish optimal decision-theoretic
properties for relative surprise inferences.

The idea of measuring surprise based on how beliefs change from \textit{a
priori} to \textit{a posteriori} and using this for inference, has arisen in
other discussions. For example, see Baldi and Itti (2010) for the use and
development of this idea in the context of learning.

While hpd and relative surprise inferences may seem quite natural, another
ingredient is often added to the formulation of a statistical problem, namely,
a loss function. For this we have an action space $\Psi,$ a function
$\Psi:\Theta\rightarrow\Psi,$ such that $\Psi(\theta)$ is the correct action
when $\theta$ is true, and a loss function $L:\Theta\times\Psi\rightarrow
\lbrack0,\infty)$ satisfying $L(\theta,\Psi(\theta))=0,$ i.e., there is no
loss when we take the correct action. The goal of a statistical decision
analysis is then to find a decision function $\delta:\mathcal{X}%
\rightarrow\Psi$ that minimizes the prior risk $r(\delta)=\int_{\Theta}%
\int_{\mathcal{X}}L(\theta,\delta(x))f_{\theta}(x)\pi(\theta)\,\mu
(dx)\,\nu(d\theta)=\int_{\mathcal{X}}r(\delta\,|\,x)m(x)\,\mu(dx)$ where
$r(\delta\,|\,x)=\int_{\Theta}L(\theta,\delta(x))\pi(\theta\,|\,x)\,\nu
(d\theta)$ is the posterior risk. Such a $\delta$ is called a Bayes rule and
clearly a $\delta$ that minimizes $r(\delta\,|\,x)$ for each $x$ is a Bayes
rule. Further discussion of decision theory can be found in Berger (1985).

As noted in Bernardo (2005) a decision formulation also leads to credible
regions for $\psi$, namely, a $\gamma$\textit{-lowest posterior loss credible
region} is defined by
\begin{equation}
L_{\gamma}(x)=\{\psi:r(\psi\,|\,x)\leq l_{\gamma}(x)\} \label{eq3}%
\end{equation}
where $l_{\gamma}(x)=\inf\{k:\int_{\{\psi_{0}:r(\psi_{0}\,|\,x)\leq k\}}%
\pi_{\Psi}(\psi\,|\,x)\,\nu_{\Psi}(d\psi)\geq\gamma\}.$ Note that $\psi$ in
(\ref{eq3}) is interpreted as the decision function that takes the value
$\psi$ constantly in $x.$ Clearly as $\gamma\rightarrow0$ the set $L_{\gamma
}(x)$ converges to the value of a Bayes rule at $x.$ For example, with
quadratic loss the Bayes rule is given by the posterior mean and a $\gamma
$-lowest posterior loss region is the smallest sphere centered at the mean
containing at least $\gamma$ of the posterior probability.

Typically, in the continuous context, Bayes rules will not be invariant under
reparameterizations. Robert (1996) recommended using the intrinsic loss
function based on a measure of distance between sampling distributions as
Bayes rules with respect to such losses are invariant. Bernardo (2005)
recommended using the intrinsic loss function based on the Kullback-Leibler
divergence $KL(f_{\theta},f_{\theta^{\prime}})$ between $f_{\theta}$ and
$f_{\theta^{\prime}}.$ When $\psi=\theta$ the intrinsic loss function is given
by $L(\theta,\theta^{\prime})=\min(KL(f_{\theta},f_{\theta^{\prime}%
}),KL(f_{\theta^{\prime}},f_{\theta})).$ For a general marginal parameter
$\psi$ the intrinsic loss function is defined by $L(\theta,\psi)=\inf
_{\theta^{\prime}\in\Psi^{-1}\{\psi\}}L(\theta,\theta^{\prime}).$

It can be shown, for example see Bernardo and Smith (2000) and Section 4, that
hpd inferences arise as the limits of Bayes rules via a sequence of loss
functions given by%
\begin{equation}
L_{\lambda}(\theta,\psi)=I(\Psi(\theta)\not \in B_{_{\lambda}}(\psi))
\label{eq4}%
\end{equation}
where $\lambda>0$ and $B_{\lambda}(\Psi(\theta))$ is the ball of radius
$\lambda$ centered at $\psi.$ As previously noted these inferences are not
invariant under reparameterizations. It is our purpose here to show that
relative surprise inferences also arise via a sequence of loss functions
similar to (\ref{eq4}) but based on the prior. So the loss functions are also
in a sense intrinsic but based on the prior and not the sampling model, as
with the intrinsic loss function.

In Section 2 we develop the prior-based loss function and show that
$\psi_{\text{LRSE}}$ is a Bayes rule when $\Psi$ is finite. In Sections 3 and
4 we extend this result to show that $\psi_{\text{LRSE}}$ is generally a limit
of Bayes rules. In Section 5 we discuss prediction problems and in Section 6
show that relative surprise regions are limits of $\gamma$-lowest\textit{
}posterior loss credible regions.

It is easy to see that the class of relative surprise credible regions
$\{C_{\gamma}(x):\gamma\in\lbrack0,1]\}$ for $\psi$ is independent of the
marginal prior $\pi_{\Psi}.$ We note, however,\ that when we specify a
$\gamma\in\lbrack0,1],$ the set $C_{\gamma}(x)$ does depend on $\pi_{\Psi}$
through $c_{\gamma}(x).$ So the form of relative surprise inferences about
$\psi$ is completely robust to the choice of $\pi_{\Psi}$ but the
quantification of the uncertainty in the inferences is not. For example, when
$\psi=\Psi(\theta)=\theta,$ then $\theta_{\text{LRSE}}(x)$ is the MLE while,
in general, $\psi_{\text{LRSE}}(x)$ is the maximizer of the integrated
likelihood where we have integrated out nuisance parameters via the
conditional prior given $\psi.$ Similarly, relative surprise regions are
likelihood regions in the case of the full parameter, and integrated
likelihood regions generally. As such, the results derived in this paper
establish that likelihood inferences are essentially Bayesian in character. We
note, however, that a relative belief ratio $\pi_{\Psi}(\psi_{0}%
\,|\,x)/\pi_{\Psi}(\psi_{0}),$ while proportional to an integrated likelihood,
has an interpretation as a change in belief and cannot be multiplied by an
arbitrary positive constant, as with a likelihood, without losing this interpretation.

In Le Cam (1953) it is shown that the MLE is asymptotically Bayes but this is
for a fixed loss function, with increasing amounts of data and a sequence of
priors. In this paper the amount of data and the prior are fixed but we may
require a sequence of loss functions, to show that the MLE is a limit of Bayes
rules. Berger, Liseo and Wolpert (1999) discuss maximum integrated likelihood
estimates where\ default or noninformative priors are used to integrate out
nuisance parameters and show good properties for this approach. Aitkin (2010)
develops an approach to assessing hypotheses using the posterior distribution
of likelihood ratios that is based on earlier work by Dempster (1973). As that
approach does not use integrated likelihoods and, as of this time, doesn't
have a decision-theoretic formulation, it is quite different than what we
discuss here.

\section{Estimation from Prior-based Loss Functions: The Finite Case}

\label{sec:prior-loss}

The following theorem presents the basic definition of the loss function when
$\Psi$ is finite and establishes an important optimality result. For more
general situations we will need to modify this loss function
slightly.\smallskip

\noindent\textbf{Theorem 1}. Suppose that $\pi_{\Psi}(\psi)>0$ for every
$\psi\in\Psi$ and that $\Psi$ is finite with $\nu_{\Psi}$ equal to counting
measure. Then for the loss function%
\begin{equation}
L(\theta,\psi)=\frac{I(\Psi(\theta)\neq\psi)}{\pi_{\Psi}(\Psi(\theta))}
\label{eq5}%
\end{equation}
a Bayes rule is given by $\psi_{\text{LRSE}}.$

\noindent Proof: We have that\newpage\
\begin{align}
r(\delta\,|\,x) &  =\int_{\Theta}\frac{I(\Psi(\theta)\neq\delta(x))}{\pi
_{\Psi}(\Psi(\theta))}\pi(\theta\,|\,x)\,\nu(d\theta)=\int_{\Psi}\frac
{I(\psi\neq\delta(x))}{\pi_{\Psi}(\psi)}\pi_{\Psi}(\psi\,|\,x)\,\nu_{\Psi
}(d\psi)\nonumber\\
&  =\int_{\Psi}\frac{\pi_{\Psi}(\psi\,|\,x)}{\pi_{\Psi}(\psi)}\,\nu_{\Psi
}(d\psi)-\frac{\pi_{\Psi}(\delta(x)\,|\,x)}{\pi_{\Psi}(\delta(x))}%
.\label{eq5a}%
\end{align}
Since $\Psi$ is finite, the first term in (\ref{eq5a}) is finite and a Bayes
rule at $x$ is given by the value $\delta(x)$ that maximizes the second
term.\ Therefore, $\psi_{\text{LRSE}}(x)$ is a Bayes rule.\smallskip

\noindent From (\ref{eq5a}) the prior risk of $\delta$ is
\begin{equation}
r(\delta)=\#(\Psi)-E_{M}(\pi_{\Psi}(\delta(x)\,|\,x)/\pi_{\Psi}(\delta
(x)))=\sum_{\psi}M_{\psi}(\delta(x)\neq\psi) \label{eq5b}%
\end{equation}
where $E_{M}$ denotes expectation with respect to the prior predictive and
$M_{\psi}$ is the probability measure on $\mathcal{X}$ obtained by averaging
$P_{\theta}$ using the conditional prior given that $\Psi(\theta)=\psi,$
namely, $M_{\psi}(A)=\int_{\Psi^{-1}\{\psi\}}P_{\theta}(A)\,\Pi(d\theta
\,|\,\Psi(\theta)=\psi).$ Therefore, finding a Bayes rule with respect to
(\ref{eq5}) is equivalent to finding $\delta$ that maximizes $E_{M}(\pi_{\Psi
}(\delta(x)\,|\,x)/\pi_{\Psi}(\delta(x))).$ So a Bayes rule maximizes the
prior expected relative belief ratio evaluated at the estimate and it is clear
that the LRSE is a Bayes rule as it maximizes the relative belief ratio for
each $x.$

If instead we take the loss function to be $I(\Psi(\theta)\neq\psi),$ then
virtually the same proof establishes that $\psi_{\text{MAP}}$ is a Bayes rule.
The prior risk for this loss function and estimator $\delta$ can be written
as
\begin{equation}
\sum_{\psi}M_{\psi}(\delta(x)\neq\psi)\pi_{\Psi}(\psi) \label{eq5c}%
\end{equation}
which is the prior probability of making an error. Both $I(\Psi(\theta
)\neq\psi)$ and (\ref{eq5}) are two-valued loss functions but, when we make an
incorrect decision, the loss is constant in $\Psi(\theta)$ for $I(\Psi
(\theta)\neq\psi)$ while it equals the reciprocal of the prior probability of
$\Psi(\theta)$ for (\ref{eq5}). So (\ref{eq5}) penalizes an incorrect decision
much more severely when the true value of $\Psi(\theta)$ is in the tails of
the prior. This makes sense as we would want to override the effect of the
prior when the prior is not placing appreciable mass at the true value. Note
that $\psi_{\text{MAP}}=\psi_{\text{LRSE}}$ when $\Pi_{\Psi}$ is uniform.

As we have already noted $\pi_{\Psi}(\psi\,|\,x)/\pi_{\Psi}(\psi)$ is
proportional to the integrated likelihood of $\psi$ when we integrate the
likelihood with respect to the conditional prior of $\theta$ given $\psi.$ So,
under the conditions of Theorem 1, we have shown that the maximum integrated
likelihood estimator is a Bayes rule. Furthermore, the Bayes rule is the same
for every choice of $\pi_{\Psi}$ and only depends on the full prior through
the conditional prior placed on the nuisance parameters. When $\psi=\theta$
then $\psi_{\text{LRSE}}(x)$ is the MLE of $\theta$ and so the MLE of $\theta$
is a Bayes rule for every prior $\pi.$

We consider an application.\smallskip\newpage

\noindent\textbf{Example 1.} \textit{Classification}

For a classification problem we have $k$ categories $\{\psi_{1},\ldots
,\psi_{k}\}$ prescribed by some function $\Psi,$ where $\pi_{\Psi}(\psi
_{i})>0$ for each $i.$ Based on observed data $x$ we want to classify the data
as having come from one of the distributions in the classes specified by
$\Psi^{-1}\{\psi_{i}\}.$

The standard Bayesian solution to this problem is to use $\psi_{\text{MAP}%
}(x)$ as the classifier. From (\ref{eq5c}) we have that $\psi_{\text{MAP}}(x)$
minimizes the prior probability of misclassification. Note that $M_{\psi
}(\delta(x)\neq\psi)$ is the prior probability of a misclassification given
that $\psi$ is the correct class and (\ref{eq5c}) is the weighted average of
these probabilities where the weights are given by the prior probabilities of
the $\psi.$ We see from (\ref{eq5b}) that $\psi_{\text{LRSE}}(x)$ is instead
minimizing the sum over $\psi$ of the probabilities of misclassification given
that $\psi$ is the correct class. So the essence of the difference between
these two approaches in this problem is that $\psi_{\text{LRSE}}(x)$ treats
the errors of misclassification equally while $\psi_{\text{MAP}}(x)$ weights
them by their prior probabilities of occurrence.

We note that (\ref{eq5b}) is an upper bound on (\ref{eq5c}). So if the Bayes
risk for loss function (\ref{eq5}) is small, the prior risk of $\psi
_{\text{LRSE}}(x),$ with respect to the loss function $I(\Psi(\theta)\neq
\psi),$ is also small, i.e., when using $\psi_{\text{LRSE}}(x)$ the overall
prior probability of a misclassification will also be small.

In general, it seems appropriate to be concerned with minimizing each of the
probabilities $M_{\psi}(\delta(x)\neq\psi)$ and not downweight those
corresponding to $\psi$ values that have small prior probability. As a
specific simple example suppose $k=2$ and $x\sim$ Binomial$(\psi_{1})$ or
$x\sim$ Binomial$(\psi_{2})$ with $\pi(\psi_{1})=1-\epsilon$ and $\pi(\psi
_{2})=\epsilon.$ After observing $x$ we want to classify the observation. For
example, $\psi_{i}$ could be the probability of a diagnostic test for a
disease indicating that the disease is present. We suppose that $\psi_{1}$ is
the probability of a positive diagnostic test for the nondiseased population
while $\psi_{2}$ is this probability for the diseased population. Further
suppose that $\psi_{1}/\psi_{2}$ is very small, indicating that the test is
successful in identifying the disease while not yielding many false positives,
and suppose $\epsilon$ is very small, indicating that the disease is very
rare. We have that $\pi(\psi_{1}\,|\,1)=\psi_{1}(1-\epsilon)/(\psi
_{1}(1-\epsilon)+\psi_{2}\epsilon)$ and $\pi(\psi_{1}\,|\,0)=(1-\psi
_{1})(1-\epsilon)/((1-\psi_{1})(1-\epsilon)+(1-\psi_{2})\epsilon).$ Therefore,
$\psi_{\text{MAP}}(1)=\psi_{1}$ if $\psi_{1}/\psi_{2}>\epsilon/(1-\epsilon)$
and is $\psi_{2}$ otherwise, while $\psi_{\text{MAP}}(0)=\psi_{1}$ if
$(1-\psi_{1})/(1-\psi_{2})>\epsilon/(1-\epsilon)$ and is $\psi_{2}$ otherwise.
Also $\psi_{\text{LRSE}}(1)=\psi_{1}$ if $\psi_{1}>\psi_{2}$ and is $\psi_{2}$
otherwise, while $\psi_{\text{LRSE}}(0)=\psi_{1}$ if $(1-\psi_{1})>(1-\psi
_{2})$ and is $\psi_{2}$ otherwise. So we see from this that $\psi
_{\text{MAP}}$ will always classify a person to the nondiseased population
when $\epsilon$ is small enough, e.g., take $\psi_{1}=0.05,\psi_{2}=0.80,$ and
$\epsilon<0.0566.$ By contrast, in this situation, $\psi_{\text{LRSE}}$ will
always classify an individual with a positive test to the diseased population
and to the nondiseased population for a negative test. Now $M_{\psi_{i}}$ is
the Binomial$(\psi_{i})$ distribution, so when $\psi_{1}<\psi_{2}$ and
$\epsilon$ is small enough
\begin{align*}
M_{\psi_{1}}(\psi_{\text{MAP}}  &  \neq\psi_{1})+M_{\psi_{2}}(\psi
_{\text{MAP}}\neq\psi_{2})=0+1=1,\\
M_{\psi_{1}}(\psi_{\text{LRSE}}  &  \neq\psi_{1})+M_{\psi_{2}}(\psi
_{\text{LRSE}}\neq\psi_{2})=\psi_{1}+(1-\psi_{2})<1.
\end{align*}
This illustrates clearly the difference between these two procedures as
$\psi_{\text{LRSE}}$ does vastly better than $\psi_{\text{MAP}}$ on the
diseased population when $\psi_{1}$ is small and $\psi_{2}$ is large as would
be the case for a good diagnostic. Of course $\psi_{\text{MAP}}$ minimizes the
overall error rate but at the price of ignoring the most important class in
this problem. Note that this example can be extended to the situation where we
need to estimate the $\psi_{i}$ based on samples from the respective
populations but this will not materially affect the overall conclusions. Also
see Example 3 where $\epsilon$ is considered unknown.\smallskip

In a general estimation problem an estimator $\delta$ is unbiased with respect
to a loss function $L$ if $E_{\theta}(L(\theta^{\prime},\delta(x)))\geq
E_{\theta}(L(\theta,\delta(x)))$ for all $\theta^{\prime},\theta\in\Theta.$
This says that on average $\delta(x)$ is closer to the true value than any
other value when we interpret $L(\theta,\delta(x))$ as a measure of distance
between the estimate and what is being estimated. A reasonable definition of
\textit{Bayesian unbiasedness} for $\delta$ with respect to \thinspace$L$ is
thus obtained by requiring that
\[
\int_{\Theta}\int_{\Theta}E_{\theta}(L(\theta^{\prime},\delta(x)))\,\Pi
(d\theta)\,\Pi(d\theta^{\prime})\geq\int_{\Theta}E_{\theta}(L(\theta
,\delta(x)))\,\Pi(d\theta)=r(\delta).
\]
Here we are thinking of $\theta^{\prime}$ as a false value generated from the
prior independently of the true value $\theta$ so $\theta^{\prime}$ has no
connection with the data. Therefore, $\delta$ is Bayesian unbiased if on
average $\delta(x)$ is closer to the true value than a false value. In Section
3 we prove that $\psi_{\text{LRSE}}$ is Bayesian unbiased with respect to a
general class of loss functions that includes both (\ref{eq5}) and
$I(\Psi(\theta)\neq\psi).$

\section{Estimation from Prior-based Loss Functions: The Countably Infinite
Case}

The loss function (\ref{eq5}) does not provide meaningful results when $\Psi$
is infinite as (\ref{eq5b}) shows that $r(\delta)$ will be infinite. So we
modify (\ref{eq5}) via a parameter $\eta>0$ and define the loss function
\begin{equation}
L_{\eta}(\theta,\psi)=\frac{I(\Psi(\theta)\neq\psi)}{\max(\eta,\pi_{\Psi}%
(\Psi(\theta)))} \label{eq6}%
\end{equation}
and note that $L_{\eta}$ is a bounded function of $(\theta,\psi).$ This loss
function is like (\ref{eq5}) but does not allow for arbitrarily large losses.
Without loss of generality we can restrict $\eta$ to a sequence of values
converging to 0. We prove the following result in the Appendix.\smallskip

\noindent\textbf{Theorem 2}. Suppose that $\pi_{\Psi}(\psi)>0$ for every
$\psi\in\Psi,$ that $\Psi$ is countable with $\nu_{\Psi}$ equal to counting
measure and that $\psi_{\text{LRSE}}(x)$ is the unique maximizer of $\pi
_{\Psi}(\psi\,|\,x)/\pi_{\Psi}(\psi)$ for all $x.$ For the loss function
(\ref{eq6}) and Bayes rule $\delta_{\eta},$ then $\delta_{\eta}(x)\rightarrow
\psi_{\text{LRSE}}(x)$ as $\eta\rightarrow0,$ for every $x\in\mathcal{X}%
.$\smallskip

\noindent The proof of Theorem also establishes the following
result.\smallskip

\noindent\textbf{Corollary 3}. For all sufficiently small $\eta$ the value of
the Bayes rule at $x$ is given by $\psi_{\text{LRSE}}(x).$\smallskip

\noindent If instead we take the loss function to be $I(\Psi(\theta)\neq
\psi),$ then virtually the same proof as in Theorem 1 establishes that
$\psi_{\text{MAP}}$ is a Bayes rule.

We now investigate the unbiasedness of $\psi_{\text{LRSE}}(x).$ For this we
consider loss functions of the form
\begin{equation}
L(\theta,\psi)=I(\Psi(\theta)\not =\psi)h(\Psi(\theta)) \label{eq6c}%
\end{equation}
for some nonnegative function $h$ which satisfies $\int_{\Theta}h(\Psi
(\theta))\,\Pi(d\theta)<\infty.$ This class of loss functions includes
(\ref{eq5}) when $\Psi$ is finite, (\ref{eq6}) and $I(\Psi(\theta)\not =%
\psi).$ We have the following result.\smallskip

\noindent\textbf{Theorem 4}. If $\Psi$ is countable, then $\psi_{\text{LRSE}%
}(x)$ is Bayesian unbiased under the loss function (\ref{eq6c}).

\noindent Proof: The prior risk of $\delta$ is given by
\begin{align*}
r(\delta)  &  =\int_{\Theta}\int_{\mathcal{X}}L(\theta,\delta(x))\,P_{\theta
}(dx)\,\Pi(d\theta)\\
&  =\int_{\Theta}\int_{\mathcal{X}}[h(\Psi(\theta))-I(\Psi(\theta
)=\delta(x))h(\Psi(\theta))]\,P_{\theta}(dx)\,\Pi(d\theta)\\
&  =\int_{\Theta}h(\Psi(\theta))\,\Pi(d\theta)-\int_{\mathcal{X}}\int_{\Theta
}I(\Psi(\theta)=\delta(x))h(\Psi(\theta))\,\Pi(d\theta\,|\,x)\,M(dx)\\
&  =\int_{\Theta}h(\Psi(\theta))\,\Pi(d\theta)-\int_{\mathcal{X}}%
h(\delta(x))\pi_{\Psi}(\delta(x)\,|\,x)\,M(dx)
\end{align*}
and%
\begin{align*}
&  \int_{\Theta}\int_{\Theta}\int_{\mathcal{X}}L(\theta^{\prime}%
,\delta(x))\,P_{\theta}(dx)\,\Pi(d\theta)\,\Pi(d\theta^{\prime})\\
&  =\int_{\Theta}\int_{\Theta}\int_{\mathcal{X}}[h(\Psi(\theta^{\prime
}))-I(\Psi(\theta^{\prime})=\delta(x))h(\Psi(\theta^{\prime}))]\,P_{\theta
}(dx)\,\Pi(d\theta)\,\Pi(d\theta^{\prime})\\
&  =\int_{\Theta}h(\Psi(\theta))\,\Pi(d\theta)-\int_{\mathcal{X}}%
h(\delta(x))\pi_{\Psi}(\delta(x))\,M(dx).
\end{align*}
Therefore, $\delta$ is Bayesian unbiased if and only if
\begin{equation}
\int_{\mathcal{X}}h(\delta(x))[\pi_{\Psi}(\delta(x)\,|\,x)-\pi_{\Psi}%
(\delta(x))]\,M(dx)\geq0. \label{eq6d}%
\end{equation}
It is a consequence of results proved in Evans and Shakhatreh (2008) that it
is always true that $\pi_{\Psi}(\psi_{\text{LRSE}}(x)\,|\,x)/\pi_{\Psi}%
(\psi_{\text{LRSE}}(x))\geq1$ and this establishes the result. This can also
be seen by noting that $\pi_{\Psi}(\cdot\,|\,x)/\pi_{\Psi}(\cdot)$ is the
density of $\Pi_{\Psi}(\cdot\,|\,x)$ with respect to $\Pi_{\Psi}$ and so we
must have that the maximum of this density is greater than or equal to
1.\smallskip

\noindent The proof gives a sufficient condition for Bayesian unbiasedness
with respect to the loss (\ref{eq6c}).\smallskip

\noindent\textbf{Corollary 5}. $\delta$ is Bayesian unbiased if $\pi_{\Psi
}(\delta(x)\,|\,x)\geq\pi_{\Psi}(\delta(x))$ for all $x.$\smallskip

At this point we have neither a proof of the Bayesian unbiasedness of
$\psi_{\text{MAP}}$ with respect to $I(\Psi(\theta)\not =\psi),$ nor a
counterexample although we suspect that it is not. We do know, however, that
$\psi_{\text{MAP}}$ is Bayesian unbiased with respect to $I(\Psi(\theta
)\not =\psi)$ whenever $\Pi_{\Psi}$ is uniform because in that case
$\psi_{\text{MAP}}=\psi_{\text{LRSE}}.$ It is also clear from (\ref{eq6c})
that $\psi_{\text{LRSE}}$ possesses a very strong property as the integrand is
always nonnegative when $\delta=\psi_{\text{LRSE}}.$ In light of this we refer
to an estimator possessing this property as being \textit{uniformly (in }%
$x$\textit{) Bayesian unbiased}.

\section{Estimation from Prior-based Loss Functions: The Continuous Case}

When $\psi$ has a continuous prior distribution the argument in Theorem 2 does
not work as $\Pi_{\Psi}(\{\delta(x)\}\,|\,x)=0.$ There are several possible
ways to proceed here but we consider a discretization of the problem that uses
Theorem 2.\ For this we will assume that the spaces involved are locally
Euclidean, mappings are sufficiently smooth and take the support measures to
be the analogs of Euclidean volume on the respective spaces. Further details
on the mathematical requirements underlying these assumptions can be found in
Tjur (1974) where spaces are taken to be Riemann manifolds. While the argument
we provide applies quite generally, we simplify this here by taking all spaces
to be open subsets of Euclidean spaces and the support measures to be
Euclidean volume on these sets.

For each $\lambda>0$ we discretize the set $\Psi$ via a countable partition
$\{B_{\lambda}(\psi):\psi\in\Psi\}$ where $\psi\in B_{\lambda}(\psi),\Pi
_{\Psi}(B_{\lambda}(\psi))>0,$ $\sup_{\psi\in\Psi}$diam$(B_{\lambda}%
(\psi))\rightarrow0$ as $\lambda\rightarrow0.$ For example, the $B_{\lambda
}(\psi)$ could be equal volume rectangles in $R^{k}.$ Further, we assume that
$\Pi_{\Psi}(B_{\lambda}(\psi))/\nu_{\Psi}(B_{\lambda}(\psi))\rightarrow
\pi_{\Psi}(\psi)$ as $\lambda\rightarrow0$ for every $\psi.$ This will hold
whenever $\pi_{\Psi}$ is continuous everywhere and $B_{\lambda}(\psi)$
converges nicely to $\{\psi\}$ as $\lambda\rightarrow0$ (see Rudin (1974),
Chapter 8 for the definition of `converges nicely'). Let $\psi_{\lambda}%
(\psi)\in B_{\lambda}(\psi)$ be such that $\psi_{\lambda}(\psi^{\prime}%
)=\psi_{\lambda}(\psi)$ whenever $\psi^{\prime}\in B_{\lambda}(\psi)$ and
$\Psi_{\lambda}=\{\psi_{\lambda}(\psi):\psi_{\lambda}(\psi)\in B_{\lambda
}(\psi)\}$ be the discretized version of $\Psi.$ Note that one point is chosen
in each $B_{\lambda}(\psi).$ We will call this a \textit{regular
discretization} of $\Psi.$ The discretized prior on $\Psi_{\lambda}$ is
$\pi_{\Psi,\lambda}(\psi_{\lambda}(\psi))=\Pi_{\Psi}(B_{\lambda}(\psi))$ and
the discretized posterior is $\pi_{\Psi,\lambda}(\psi_{\lambda}(\psi
)\,|\,x)=\Pi_{\Psi}(B_{\lambda}(\psi)\,|\,x).$

We define the loss function for the discretized problem just as for Theorem 2,
by%
\begin{equation}
L_{\lambda,\eta}(\theta,\psi_{\lambda}(\psi))=\frac{I(\psi_{\lambda}%
(\Psi(\theta))\neq\psi_{\lambda}(\psi))}{\max(\eta,\pi_{\Psi,\lambda}%
(\psi_{\lambda}(\Psi(\theta))))} \label{eq8}%
\end{equation}
and denote a Bayes rule for this problem by $\delta_{\lambda,\eta}(x).$ In
this case we not only need that $\psi_{\text{LRSE}}(x)$ is the unique
maximizer of $\pi_{\Psi}(\psi\,|\,x)/\pi_{\Psi}(\psi),$ but we cannot allow
$\pi_{\Psi}(\psi\,|\,x)/\pi_{\Psi}(\psi)$ to come arbitrarily close to its
maximum outside a neighborhood of $\psi_{\text{LRSE}}(x).$ It is clear that
when this does not hold then we are in a pathological situation that will not
apply in a typical application. The following result is proved in the
Appendix.\smallskip

\noindent\textbf{Theorem 6}. Suppose that $\pi_{\Psi}$ is positive and
continuous and we have a regular discretization of $\Psi.$ Further suppose
that $\psi_{\text{LRSE}}(x)$ is the unique maximizer of $\pi_{\Psi}%
(\psi\,|\,x)/\pi_{\Psi}(\psi)$ and for any $\epsilon>0$\newpage\
\[
\sup_{\{\psi:||\psi-\psi_{\text{LRSE}}(x)||\geq\epsilon\}}\frac{\pi_{\Psi
}(\psi\,|\,x)}{\pi_{\Psi}(\psi)}<\frac{\pi_{\Psi}(\psi_{\text{LRSE}%
}(x)\,|\,x)}{\pi_{\Psi}(\psi_{\text{LRSE}}(x))}.
\]
Then, there exists $\eta(\lambda)>0$ such that a Bayes rule $\delta
_{\lambda,\eta(\lambda)}(x)$ converges to $\psi_{\text{LRSE}}(x)$ as
$\lambda\rightarrow0$ for all $x.$\smallskip

\noindent Theorem 6 says that $\psi_{\text{LRSE}}$ is a limit of Bayes rules.
So when $\Psi(\theta)=\theta$ we have the result that the MLE is a limit of
Bayes rules and more generally the maximum integrated likelihood estimator is
a limit of Bayes rules.

Now let $\hat{\psi}_{\lambda}(x)$ be the LRSE of the discretized problem,
i.e., $\hat{\psi}_{\lambda}(x)$ maximizes $\Pi_{\Psi}(B_{\lambda}%
(\psi)\,|\,x)/\Pi_{\Psi}(B_{\lambda}(\psi))$ as a function of $\psi\in
\Psi_{\lambda}.$ The following result is proved in the Appendix.\smallskip

\noindent\textbf{Corollary 7}. $\hat{\psi}_{\lambda}$ converges to
$\psi_{\text{LRSE}}$ as $\lambda\rightarrow0.$\smallskip

\noindent Note that by Theorem 4, $\hat{\psi}_{\lambda}$ is uniformly Bayesian
unbiased for the discretized problem. Therefore, $\psi_{\text{LRSE}}$ is the
limit of uniformly Bayesian unbiased estimators.

By similar arguments we can establish an analog of Theorem 6 for
$\psi_{\text{MAP}}$ using the loss function given by (\ref{eq4}). Actually in
this case a simpler development can be followed in certain situations. For
this note that the posterior risk of $\delta$ is given by $1-\Pi_{\Psi
}(B_{\lambda}(\delta(x))\,|\,x)=1-\pi_{\Psi}(\delta^{\prime}(x)\,|\,x)\nu
_{\Psi}(B_{\lambda}(\delta(x)))$ for some $\delta^{\prime}(x)\in B_{\lambda
}(\delta(x)).$ Now suppose we take $B_{\lambda}(\psi)$ to be a sphere of
radius $\lambda$ centered at $\psi.$ Suppose further that for each
$\epsilon>0$ there exists a $\lambda(\epsilon)>0$ such that when $||\psi
-\psi_{\text{MAP}}(x)||>\lambda(\epsilon)$ then $\pi_{\Psi}(\psi
\,|\,x)<\inf_{\psi^{\prime}\in B_{\lambda(\epsilon)}(\psi_{\text{MAP}}(x))}%
\pi_{\Psi}(\psi^{\prime}\,|\,x).$ Since $\nu_{\Psi}(B_{\lambda}(\psi))$ is
constant we have that a Bayes rule $\delta_{\lambda}$ must then satisfy
$||\delta_{\lambda}(x)-\psi_{\text{MAP}}(x)||<\epsilon$. So we have proved
that $\psi_{\text{MAP}}$ is a limit of Bayes rules. By contrast, for the loss
function $I(\Psi(\theta)\not \in B_{_{\lambda}}(\psi))/\Pi_{\Psi}(B_{\lambda
}(\Psi(\theta)))$ the posterior risk of $\delta$ is given by $\int_{\Psi}%
\{\Pi_{\Psi}(B_{\lambda}(\psi))\}^{-1}\Pi_{\Psi}(d\psi\,|\,x)-\int%
_{B_{\lambda}(\delta(x))}\{\Pi_{\Psi}(B_{\lambda}(\psi))\}^{-1}\Pi_{\Psi
}(d\psi\,|\,x).$ The simpler approach is not available in this case because
the first term is unbounded.

We consider now an important example.\smallskip

\noindent\textbf{Example 2.} \textit{Regression (estimation)}

Suppose that we have $y=X\beta+e$ where $y\in R^{n},X\in R^{n\times k}$ is
fixed, $\beta\in R^{n\times k},$ and $e\sim N_{n}(0,\sigma^{2}I).$ We will
assume that $\sigma^{2}$ is known to simplify the discussion. Let $\pi$ be a
prior density for $\beta.$ Then having observed $(X,y),$ $\beta_{\text{LRSE}%
}(x)=b=(X^{\prime}X)^{-1}X^{\prime}y$ which is the MLE of $\beta.$ It is
interesting to contrast this result with what might be considered more
standard Bayesian estimates such as the posterior mode or posterior mean. For
example, suppose that $\beta\sim N_{k}(0,\tau^{2}I).$ Then the posterior
distribution of $\beta\ $is $N_{k}(\mu_{\text{post}}(\beta),\Sigma
_{\text{post}}(\beta))$ where
\[
\mu_{\text{post}}(\beta)=\Sigma_{\text{post}}(\beta)\sigma^{-2}X^{\prime
}Xb,\text{ }\Sigma_{\text{post}}(\beta)=(\tau^{-2}I+\sigma^{-2}X^{\prime
}X)^{-1}%
\]
and the posterior mean and modal estimates of $\beta$ are both equal to
$\mu_{\text{post}}(\beta).$ Writing the spectral decomposition of $X^{\prime
}X$ as $X^{\prime}X=Q\Lambda Q^{\prime}$ we have that
\[
||\mu_{\text{post}}(\beta)||=||(I+(\sigma^{2}/\tau^{2})\Lambda^{-1}%
)^{-1}Q^{\prime}b||.
\]
Since $||b||=||Q^{\prime}b||$ and $1/(1+\sigma^{2}/(\tau^{2}\lambda_{i}))<1$
for each $i,$ we see that $\mu_{\text{post}}(\beta)$ moves the MLE towards the
prior mean 0. This is often cited as a positive attribute of these estimates
but consider the situation where the true value of $\beta$ lies in the tails
of the prior. In that case it is certainly wrong to move $\beta$ towards the
prior mean. When $\tau^{2}$ is chosen very large, so we avoid the possibility
that the true value of $\beta$ lies in the tails of the prior, then the MLE
and the posterior mean are virtually the same. It makes sense to choose
$\tau^{2}>\sigma^{2}$ as this says we have less prior information about a
$\beta_{i}$ than the amount we learn about $\beta$ from a single observation.
So it is not clear that shrinking the MLE is necessarily a good thing
particularly as this requires giving up invariance.

Suppose now we want to estimate $\psi=w^{\prime}\beta$ for some setting $w$ of
the predictors. The prior distribution of $\psi$ is $N(0,\sigma_{\text{prior}%
}^{2}(\psi))=N(0,\tau^{2}w^{\prime}w)$ and the posterior distribution is
$N(\mu_{\text{post}}(\psi),\sigma_{\text{post}}^{2}(\psi))=N(w^{\prime}%
\mu_{\text{post}}(\beta),x^{\prime}\Sigma_{\text{post}}(\beta)x).$ Note that
$\sigma_{\text{prior}}^{2}(\psi)-\sigma_{\text{post}}^{2}(\psi)=w^{\prime
}(\tau^{2}I-\Sigma_{\text{post}}(\beta))w=\tau^{2}w^{\prime}Q^{\prime
}(I-(I+(\tau^{2}/\sigma^{2})\Lambda)^{-1})Qw>0$ and so maximizing the ratio of
the posterior to prior densities leads to
\begin{equation}
\psi_{\text{LRSE}}(y)=(1-\sigma_{\text{post}}^{2}(\psi)/\sigma_{\text{prior}%
}^{2}(\psi))^{-1}\mu_{\text{post}}(\psi). \label{est1}%
\end{equation}
Since $\sigma_{\text{prior}}^{2}(\psi)>\sigma_{\text{post}}^{2}(\psi)$ we have
$|\psi_{\text{LRSE}}(y)|>|\mu_{\text{post}}(\psi)|$ and $\mu_{\text{post}%
}(\psi)=\psi_{\text{MAP}}(y).$ Note that when $\sigma_{\text{post}}^{2}(\psi)$
is much smaller than $\sigma_{\text{prior}}^{2}(\psi),$ in other words the
posterior is densely concentrated about $\mu_{\text{post}}(\psi),$ then
$w_{\text{LRSE}}(y)$ and $w_{\text{MAP}}(y)$ are very similar. In general
$\psi_{\text{LRSE}}(y)$ is not equal to $w^{\prime}b,$ the plug-in MLE of
$\psi,$ although $\psi_{\text{LRSE}}(y)\rightarrow w^{\prime}b$ as $\tau
^{2}\rightarrow\infty.$

\section{Prediction from Prior-based Loss Functions}

Suppose after observing $x$ we want to predict a future value $y\in
\mathcal{Y}$ where $y$ has model given by $g_{\eta(\theta)}(y\,|\,x)$ with
respect to support measure $\mu_{\mathcal{Y}}$ on $\mathcal{Y}.$ We allow for
the possibility here that the distribution of $y$ depends on $x$ and also that
$\theta$ may not index these distributions. Then we have that the joint
density of $(\theta,x,y)$ is given by $\pi(\theta)f_{\theta}(x)g_{\eta
(\theta)}(y\,|\,x)$ and after observing $x$ the conditional density of $y$ is
given by the posterior predictive density $q(y\,|\,x)=\int_{\Theta}\pi
(\theta\,|\,x)g_{\eta(\theta)}(y\,|\,x)\,\nu(d\theta)$ while the prior
predictive density of $y$ is given by $q(y)=\int_{\Theta}\int_{\mathcal{X}}%
\pi(\theta)f_{\theta}(x)g_{\eta(\theta)}(y\,|\,x)\,\mu(dx)\,\nu(d\theta).$
Therefore, the relative belief in a future value $y$ is given by
$q(y\,|\,x)/q(y)$ and we denote the maximizer of this by $y_{\text{LRSE}}(x).$

Again the LRSE arises from loss function considerations. For example, when
$\mathcal{Y}$ is finite we consider the loss function
\[
L(y,y^{\prime})=\frac{I(y\neq y^{\prime})}{q(y)},
\]
where we think of $y$ as some true value of $y$ that is concealed from us by
the future, or some other mechanism, and which we want to predict. Then the
posterior risk of a predictor $\delta:\mathcal{X}\rightarrow\mathcal{Y}$ is
given by
\[
r(\delta\,|\,x)=\int_{\mathcal{Y}}\frac{q(y\,|\,x)}{q(y)}\,\mu_{\mathcal{Y}%
}(dy)-\frac{q(\delta(x)\,|\,x)}{q(\delta(x))}%
\]
and we see that $y_{\text{LRSE}}$ is a Bayes rule. Also, the prior risk of
predictor $\delta$ is given by $r(\delta)=\sum_{y}M_{y}(\delta(x)\neq y)$
where $M_{y}$ is the conditional prior predictive of $x$ given $y$ and so
$r(\delta)$ is the sum of the conditional prediction errors given $y.$ We can
also develop results similar to Theorems 2 and 6 for the situation where
$\mathcal{Y}$ is not finite to show that $y_{\text{LRSE}}$ is a limit of Bayes rules.

We consider some examples.\smallskip

\noindent\textbf{Example 3. }\textit{Classification (prediction)}

Consider now a situation where $(x,c)$ is such that $x\,|\,c\sim f_{c}$ with
$c\sim$ Bernoulli$(\epsilon)$ where $f_{0}$ and $f_{1}$ are known (or
accurately estimated based on large samples) but $\epsilon$ is unknown with
prior $\pi.$ This is a generalization of Example 1 where $\epsilon$ was
assumed to be known. Then based on a sample $(x_{1},c_{1}),\ldots,(x_{n}%
,c_{n})$ from the joint distribution we want to predict the value $c_{n+1}$
for a newly observed $x_{n+1}.$ Therefore, $q(c)=\int_{0}^{1}(1-\epsilon
)^{1-c}\epsilon^{c}\pi(\epsilon)\,d\epsilon\ $and, if $\epsilon\sim$
Beta$(\alpha,\beta),$ the prior predictive of $c_{n+1}$ is Bernoulli$(\alpha
/(\alpha+\beta)).$ For $c_{n+1}$ the posterior predictive density is
$q(c\,|\,(x_{1},c_{1}),\ldots,(x_{n},c_{n}),x_{n+1})\propto$ $(f_{0}%
(x_{n+1}))^{1-c}(f_{1}(x_{n+1}))^{c}\int_{0}^{1}\epsilon^{n\bar{c}%
+c}(1-\epsilon)^{n(1-\bar{c})+(1-c)}\pi(\epsilon)\,d\epsilon$ with $\bar
{c}=n^{-1}\sum_{i=1}^{n}\newline c_{i}.$ With a Beta$(\alpha,\beta)$ prior for
$\epsilon,$ we have that $q(c\,|\,(x_{1},c_{1}),\ldots,(x_{1},c_{1}),x_{n+1}%
)$\newline$\propto f_{c}(x_{n+1})\Gamma\left(  \alpha+n\bar{c}+c\right)
\Gamma(\beta+n(1-\bar{c})+1-c).$ From this we see immediately that%
\begin{align}
c_{\text{MAP}}  &  =\left\{
\begin{array}
[c]{cl}%
1 & \text{if }\frac{f_{1}(x_{n+1})}{f_{0}(x_{n+1})}\frac{\left(  \alpha
+n\bar{c}\right)  }{\left(  \beta+n(1-\bar{c})\right)  }\geq1\\
0 & \text{otherwise,}%
\end{array}
\right. \nonumber\\
c_{\text{LRSE}}  &  =\left\{
\begin{array}
[c]{cl}%
1 & \text{if }\frac{f_{1}(x_{n+1})}{f_{0}(x_{n+1})}\frac{\beta\left(
\alpha+n\bar{c}\right)  }{\alpha\left(  \beta+n(1-\bar{c})\right)  }\geq1\\
0 & \text{otherwise.}%
\end{array}
\right.  \label{predclass}%
\end{align}
Note that $c_{\text{MAP}}$ and $c_{\text{LRSE}}$ are identical whenever
$\alpha=\beta.$

We can see from these formulas that a substantial difference will arise
between $c_{\text{MAP}}$ and $c_{\text{LRSE}}$ when one of $\alpha$ or $\beta$
is much bigger than the other. As in Example 1 these correspond to situations
where we believe that $\epsilon$ or $1-\epsilon$ is very small. Suppose we
take $\alpha=1$ and let $\beta$ be relatively large, as this corresponds to
knowing \textit{a priori} that $\epsilon$ is very small. Then (\ref{predclass}%
) implies that $c_{\text{MAP}}\leq c_{\text{LRSE}}$ and so $c_{\text{LRSE}}=1$
whenever $c_{\text{MAP}}=1.$ A similar conclusion arises when we take
$\beta=1$ and $\alpha<1.$

To see what kind of improvement is possible we consider a simulation. Here we
take $f_{0}$ to be a $N(0,1)$ density, $f_{1}$ to be a $N(\mu,1)$ density, let
$n=10$ and the prior on $\epsilon$ be Beta$(1,\beta).$ Table 1 presents the
Bayes risks for $c_{\text{MAP}}\ $and $c_{\text{LRSE}}$ for various choices of
$\beta$ when $\mu=1.$ When $\beta=1$ they are equivalent but we see that as
$\beta$ rises the performance of $c_{\text{MAP}}$ deteriorates while
$c_{\text{LRSE}}$ improves. Large values of $\beta$ correspond to having
information that $\epsilon$ is small. When $\beta=14$ about 0.50 of the prior
probability is to the left of 0.05, with $\beta=32$ about 0.80 of the prior
probability is to the left of 0.05, and with $\beta=100$ about 0.99 of the
prior probability is to the left of 0.05. We see that the misclassification
rates for the small group $(c=1)$ stay about the same for $c_{\text{LRSE}}$ as
$\beta$ increases while they deteriorate markedly for $c_{\text{MAP}}$ as the
MAP\ procedure basically ignores the small group.%

%TCIMACRO{\TeXButton{B}{\begin{table}[tbp] \centering}}%
%BeginExpansion
\begin{table}[tbp] \centering
%EndExpansion%
\begin{tabular}
[c]{|r|r|r|}\hline
$\beta$ & $M_{0}(c_{\text{MAP}}\neq0)+M_{1}(c_{\text{MAP}}\neq1)$ &
$M_{0}(c_{\text{LRSE}}\neq0)+M_{1}(c_{\text{LRSE}}\neq1)$\\\hline
$1$ & \multicolumn{1}{|c|}{$0.386+0.390=0.776$} &
\multicolumn{1}{|c|}{$0.386+0.390=0.776$}\\
$14$ & \multicolumn{1}{|c|}{$0.002+0.975=0.977$} &
\multicolumn{1}{|c|}{$0.285+0.380=0.665$}\\
$32$ & \multicolumn{1}{|c|}{$0.000+0.997=0.997$} &
\multicolumn{1}{|c|}{$0.292+0.349=0.641$}\\
$100$ & \multicolumn{1}{|c|}{$0.000+1.000=1.000$} &
\multicolumn{1}{|c|}{$0.300+0.324=0.624$}\\\hline
\end{tabular}
\caption{Conditional prior probabilities of misclassification for MAP and LRSE for various values of $\beta$ in Example 3 when $\alpha=1$, $\mu=1$, and $n$=10. }\label{TableKey}%
%TCIMACRO{\TeXButton{E}{\end{table}}}%
%BeginExpansion
\end{table}%
%EndExpansion

We also investigated other choices for $n$ and $\mu.$ There is very little
change as $n$ increases. When $\mu$ moves towards 0 the error rates go up and
go down as $\mu$ moves away from 0, as one would expect. Of course,
$c_{\text{LRSE}}$ always dominates $c_{\text{MAP}}.$\smallskip

\noindent\textbf{Example 4.} \textit{Regression (prediction)}

Consider the situation of \ Example 2 and suppose we want to predict a
response $z$ at the predictor value $w\in R^{k}.$ When $\beta\sim N_{k}%
(0,\tau^{2}I)$ the prior distribution of $z$ is $z\sim N(0,\sigma^{2}+\tau
^{2}w^{\prime}w)=N(0,\sigma_{\text{prior}}^{2}(z))$ and the posterior
distribution is $N(\mu_{\text{post}}(z),\sigma_{\text{post}}^{2}(z))$ where
\[
\mu_{\text{post}}(z)=w^{\prime}\mu_{\text{post}}(\beta),\text{ }%
\sigma_{\text{post}}^{2}(z)=\sigma^{2}+w^{\prime}\Sigma_{\text{post}}%
(\beta)w.
\]
To obtain $z_{\text{LRSE}}(y)$ we need to maximize the ratio of the posterior
to the prior density of $z$ and an easy calculation shows that this leads to%
\begin{equation}
z_{\text{LRSE}}(y)=(1-\sigma_{\text{post}}^{2}(z)/\sigma_{\text{prior}}%
^{2}(z))^{-1}\mu_{\text{post}}(z). \label{pred1}%
\end{equation}
Note that $\sigma_{\text{prior}}^{2}(z)-\sigma_{\text{post}}^{2}%
(z)=\sigma_{\text{prior}}^{2}(w^{\prime}\beta)-\sigma_{\text{post}}%
^{2}(w^{\prime}\beta)>0$ and so $|z_{\text{LRSE}}(y)|$\newline$>|\mu
_{\text{post}}(z)|$ and the LRSE is further from the prior mean than
$z_{\text{MAP}}(y)=\mu_{\text{post}}(z).$ Also, we see that, when
$\sigma_{\text{post}}^{2}(z)$ is small then $z_{\text{LRSE}}(y)$ and
$z_{\text{MAP}}(y)$ are very similar. Finally, comparing (\ref{est1}) and
(\ref{pred1}) we have that
\[
z_{\text{LRSE}}(y)=(\sigma_{\text{prior}}^{2}(z)/\sigma_{\text{post}}^{2}%
(\psi))w^{\prime}\psi_{\text{LRSE}}(y)=(1+\sigma^{2}/\tau^{2})\psi
_{\text{LRSE}}(y)
\]
and so the LRSE predictor at $x$ is more dispersed than the LRSE estimator of
the mean at $w$ and this makes good sense as we have to take into account the
additional variation due to prediction. By contrast $w_{\text{MAP}}%
(y)=\psi_{\text{MAP}}(y).$

\section{Regions from Prior-based Loss Functions}

We now consider the lowest posterior loss $\gamma$-credible regions that arise
from the prior-based loss functions we have considered. Let $C_{\gamma}(x)$
denote a\ $\gamma$-relative surprise region for $\psi.$ Consider first the
case where $\Psi$ is finite. We have the following result.\smallskip

\noindent\textbf{Theorem 8}. Suppose that $\pi_{\Psi}(\psi)>0$ for every
$\psi\in\Psi$ and that $\Psi$ is finite with $\nu_{\Psi}$ equal to counting
measure. Then for the loss function given by (\ref{eq5}), $C_{\gamma}(x)$ is a
$\gamma$-lowest posterior loss credible region.

\noindent Proof: From (\ref{eq3}) and (\ref{eq5a}) the $\gamma$-lowest
posterior loss credible region is%
\[
L_{\gamma}(x)=\left\{  \psi:\frac{\pi_{\Psi}(\psi\,|\,x)}{\pi_{\Psi}(\psi
)}\geq\int_{\Psi}\frac{\pi_{\Psi}(\zeta\,|\,x)}{\pi_{\Psi}(\zeta)}\,\nu_{\Psi
}(d\zeta)-l_{\gamma}(x)\right\}
\]
and $l_{\gamma}(x)=\inf\{k:\Pi_{\Psi}(\{\psi:r(\psi\,|\,x)\leq k\}\geq
\gamma\}.$ As $\int_{\Psi}(\pi_{\Psi}(z\,|\,x)/\pi_{\Psi}(z))\nu_{\Psi}(dz)$
is independent of $\psi$ it is clearly equivalent to define this region via
$C_{\gamma}(x)=\left\{  \psi:\pi_{\Psi}(\psi\,|\,x)/\pi_{\Psi}(\psi)\geq
c_{\gamma}(x)\right\}  ,$ namely, $L_{\gamma}(x)=C_{\gamma}(x).$\smallskip

Now consider the case where $\Psi$ is countable and we use loss function
(\ref{eq6}). Following the proof of Theorem 8 we see that a $\gamma$-lowest
posterior loss region takes the form%
\[
L_{\eta,\gamma}(x)=\left\{  \psi:\pi_{\Psi}(\psi\,|\,x)/\max(\eta,\pi_{\Psi
}(\psi))\geq l_{\eta,\gamma}(x)\right\}
\]
where $l_{\eta,\gamma}(x)=\sup\{k:\Pi_{\Psi}(\{\psi:\pi_{\Psi}(\psi
\,|\,x)/\max(\eta,\pi_{\Psi}(\psi))\geq k\}\geq\gamma\}.$ We prove the
following result in the Appendix.\smallskip

\noindent\textbf{Theorem 9}. Suppose that $\pi_{\Psi}(\psi)>0$ for every
$\psi\in\Psi,$ that $\Psi$ is countable with $\nu_{\Psi}$ equal to counting
measure. For the loss function (\ref{eq6}), we have that $C_{\gamma}%
(x)\subset\lim\inf_{\eta\rightarrow0}L_{\eta,\gamma}(x)$ whenever $\gamma$ is
such that $\Pi_{\Psi}(C_{\gamma}(x)\,|\,x)=\gamma$ and $\lim\sup
_{\eta\rightarrow0}L_{\eta,\gamma}(x)\subset C_{\gamma^{\prime}}(x)$ whenever
$\gamma^{\prime}>\gamma$ and $\Pi_{\Psi}(C_{\gamma^{\prime}}(x)\,|\,x)=\gamma
^{\prime}.$\smallskip

\noindent While Theorem 9 does not establish the exact convergence $\lim
_{\eta\rightarrow0}L_{\eta,\gamma}(x)=C_{\gamma}(x)$ we suspect, however, that
this does hold under quite general circumstances due to the discreteness.
Theorem 9 does show that limit points of the class of sets $L_{\eta,\gamma
}(x)$ always contain $C_{\gamma}(x)$ and their posterior probability content
differs from $\gamma$ by at most $\gamma^{\prime}-\gamma$ where $\gamma
^{\prime}>\gamma$ is the next largest value for which we have exact content.

We now consider the continuous case and suppose we have a regular
discretization. For $S^{\ast}\subset\Psi_{\lambda}=\{\psi_{\lambda}(\psi
):\psi_{\lambda}(\psi)\in B_{\lambda}(\psi)\},$ namely, $S^{\ast}$ is a subset
of a discretized version of $\Psi,$ we define the \textit{undiscretized}
version of $S^{\ast}$ to be $S=\cup_{\psi\in S^{\ast}}B_{\lambda}(\psi).$ Now
let $C_{\lambda,\gamma}^{\ast}(x)$ be the $\gamma$-relative surprise region
for the discretized problem and let $C_{\lambda,\gamma}(x)$ be its
undiscretized version. Note that in a continuous context we will consider two
sets as equal if they differ only by a set of measure 0 with respect to
$\Pi_{\Psi}.$ In the Appendix we prove the following which says that a
$\gamma$-relative surprise region for the discretized problem (after
undiscretizing) converges to the $\gamma$-relative surprise region for the
original problem.\smallskip

\noindent\textbf{Theorem 10}. Suppose that $\pi_{\Psi}$ is positive and
continuous, we have a regular discretization of $\Psi$ and $\pi_{\Psi}%
(\psi\,|\,x)/\pi_{\Psi}(\psi)$ has a continuous posterior distribution. Then
$\lim_{\lambda\rightarrow0}C_{\lambda,\gamma}(x)=C_{\gamma}(x).$\smallskip

\noindent While Theorem 10 has interest in its own right, we can use it to
prove that relative surprise regions are limits of lowest posterior loss regions.

Let $L_{\eta,\lambda,\gamma}^{\ast}(x)$ be the $\gamma$-lowest posterior loss
region obtained for the discretized problem using loss function (\ref{eq8})
and let $L_{\eta,\lambda,\gamma}(x)$ be the undiscretized version. We prove
the following result in the Appendix.\smallskip

\noindent\textbf{Theorem 11}. Suppose that $\pi_{\Psi}$ is positive and
continuous, we have a regular discretization of $\Psi$ and $\pi_{\Psi}%
(\psi\,|\,x)/\pi_{\Psi}(\psi)$ has a continuous posterior distribution. Then
$C_{\gamma}(x)=\underset{\lambda\rightarrow0}{\lim}\underset{\eta
\rightarrow0}{\lim\inf\,}L_{\eta,\lambda,\gamma}(x)=\underset{\lambda
\rightarrow0}{\lim}\underset{\eta\rightarrow0}{\lim\sup}\,L_{\eta
,\lambda,\gamma}(x).$

In Evans, Guttman, and Swartz (2006) and Evans and Shakhatreh (2008)
additional properties of relative surprise regions are developed. For example,
it is proved that a $\gamma$-relative surprise region $C_{\gamma}(x)$ for
$\psi$ satisfying $\Pi_{\Psi}(C_{\gamma}(x)\,|\,x)=\gamma$ minimizes
$\Pi_{\Psi}(B)$ among all (measurable) subsets of $\Psi$ satisfying $\Pi
_{\Psi}(B\,|\,x)\geq\gamma.$ So a $\gamma$-relative surprise region is
smallest among all $\gamma$-credible regions for $\psi$ where size is measured
using the prior measure. This property has several consequences. For example,
the prior probability that a region $B(x)\subset\Psi$ contains a false value
from the prior is given by $\int_{\Theta}\int_{\Psi}P_{\theta}(\psi\in
B(x))\,\Pi_{\Psi}(d\psi)\,\Pi(d\theta)$ where a false value is a value of
$\psi\sim\Pi_{\Psi}$ generated independently of $(\theta,x)\sim\Pi_{\Psi
}\times P_{\theta}.$ It can be proved that a $\gamma$-relative surprise region
minimizes this probability among all $\gamma$-credible regions for $\psi$ and
is always unbiased in the sense that the probability of covering a false value
is bounded above by $\gamma.$ Furthermore, a $\gamma$-relative surprise region
maximizes the relative belief ratio $\Pi_{\Psi}(B\,|\,x)/\Pi_{\Psi}(B)$ and
the Bayes factor $\Pi_{\Psi}(B\,|\,x)\Pi_{\Psi}(B^{c})/\Pi_{\Psi}%
(B^{c}\,|\,x)\Pi_{\Psi}(B)$ among all regions $B\subset\Psi$ with $\Pi_{\Psi
}(B)=\Pi_{\Psi}(C_{\gamma}(x)\,|\,x).$

While the results in this section have been concerned with obtaining credible
regions for parameters, similar results can be proved for the construction of
prediction regions.

\section{Conclusions}

Relative surprise inferences are closely related to likelihood inferences.
This together with their invariance and optimality properties make these prime
candidates as appropriate inferences in Bayesian contexts.\ This paper has
shown that relative surprise inferences arise naturally in a
decision-theoretic formulation using loss functions based on the prior. As of
yet these inferences are not typically used while MAP-based inferences, which
seem to possess few strong properties, are commonly recommended. Based on the
properties we have discussed in this paper we conclude that improvements in
inferences can be accomplished by adopting relative surprise inferences. While
we have required proper priors in this paper, limiting relative surprise
inferences, as priors become increasingly diffuse, can also be obtained and
have been discussed in the references.

Relative surprise estimation of the parameter $\psi$ is based on the relative
belief ratio $\pi_{\Psi}(\psi\,|\,x)/\pi_{\Psi}(\psi).$ As this ratio is
independent of the choice of $\pi_{\Psi},$ estimation of $\psi$ is to a
certain extent robust to the choice of prior. The role of the marginal prior
$\pi_{\Psi}$ arises in quantifying the uncertainty about the estimate of
$\psi$ through the regions $C_{\gamma}.$ So the conditional prior given
$\psi,$ together with the model and data, are used to determine the form of
any inferences about $\psi,$ while the marginal prior for $\psi,$ together
with the model and data, are used to quantify the uncertainty in these inferences.

By contrast predictions are based on the relative belief ratio
$q(y\,|\,x)/q(y)$ which is generally dependent on the full prior $\pi.$ So in
a sense predictions are less robust to the prior than estimation. On the other
hand Bayesian inferences are often advocated due to the regularizing effect of
the prior. While the relative surprise approach does not fully incorporate
such an effect for parameter estimates, the full effect is available for prediction.

\section*{Appendix}

\noindent\textbf{Proof of Theorem} 2: We have that
\begin{align}
r_{\eta}(\delta\,|\,x)  &  =\int_{\Psi}\frac{I(\psi\neq\delta(x))}{\max
(\eta,\pi_{\Psi}(\psi))}\pi_{\Psi}(\psi\,|\,x)\,\nu_{\Psi}(d\psi)\nonumber\\
&  =\int_{\Psi}\frac{\pi_{\Psi}(\psi\,|\,x)}{\max(\eta,\pi_{\Psi}(\psi))}%
\,\nu_{\Psi}(d\psi)-\frac{\pi_{\Psi}(\delta(x)\,|\,x)}{\max(\eta,\pi_{\Psi
}(\delta(x)))}. \label{eq6a}%
\end{align}
The first term in (\ref{eq6a}) is constant in $\delta(x)$ and bounded above by
$1/\eta,$ so the value of a Bayes rule at $x$ is obtained by finding
$\delta(x)$ that maximizes the second term.

Consider $\eta$ as fixed and note that
\begin{equation}
\frac{\pi_{\Psi}(\delta(x)\,|\,x)}{\max(\eta,\pi_{\Psi}(\delta(x)))}=\left\{
\begin{array}
[c]{cl}%
\frac{\pi_{\Psi}(\delta(x)\,|\,x)}{\eta} & \text{if }\eta>\pi_{\Psi}%
(\delta(x))\\
\frac{\pi_{\Psi}(\delta(x)\,|\,x)}{\pi_{\Psi}(\delta(x))} & \text{if }\eta
\leq\pi_{\Psi}(\delta(x)).
\end{array}
\right.  \label{eq6b}%
\end{equation}
There are at most finitely many values of $\psi$ satisfying $\eta\leq\pi
_{\Psi}(\psi)$ and so $\pi_{\Psi}(\psi\,|\,x)/\pi_{\Psi}(\psi)$ assumes a
maximum on this set, say at $\psi_{\eta}(x)$. There are infinitely many values
of $\psi$ satisfying $\eta>\pi_{\Psi}(\psi)$ but clearly we can find
$\eta^{\prime}<\eta$ so that $\{\psi:\eta^{\prime}<\pi_{\Psi}(\psi)<\eta\}$ is
nonempty and finite. Thus, $\pi_{\Psi}(\psi\,|\,x)$ assumes its maximum on the
set $\{\psi:\pi_{\Psi}(\psi)<\eta\}$ in the subset $\{\psi:\eta^{\prime}%
<\pi_{\Psi}(\psi)<\eta\},$ say at $\psi_{\eta}^{\prime}(x).$ Therefore, a
Bayes rule $\delta_{\eta}(x)$ is given by $\delta_{\eta}(x)=\psi_{\eta}(x)$
when $\pi_{\Psi}(\psi_{\eta}(x)\,|\,x)/\pi_{\Psi}(\psi_{\eta}(x))\geq\pi
_{\Psi}(\psi_{\eta}^{\prime}(x)\,|\,x)/\eta$ and $\delta_{\eta}(x)=\psi_{\eta
}^{\prime}(x)$ otherwise.

If $\eta>\pi_{\Psi}(\delta(x)),$ then
\[
\pi_{\Psi}(\delta(x)\,|\,x)/\eta<\pi_{\Psi}(\delta(x)\,|\,x)/\pi_{\Psi}%
(\delta(x))\leq\pi_{\Psi}(\psi_{\text{LRSE}}(x)\,|\,x)/\pi_{\Psi}%
(\psi_{\text{LRSE}}(x)).
\]
Therefore, whenever $\eta\leq\pi_{\Psi}(\psi_{\text{LRSE}}(x))$ the maximizer
of (\ref{eq6b}) is given by $\delta(x)=\psi_{\text{LRSE}}(x)$ and the result
is proved.\medskip

\noindent\textbf{Proof of Theorem 6}: Just as in Theorem 2 a Bayes rule
$\delta_{\lambda,\eta}(x)$ maximizes $\pi_{\Psi,\lambda}(\delta(x)\,|\,x)/\max
(\eta,\pi_{\Psi,\lambda}(\delta(x)))$ for $\delta(x)\in\Psi_{\lambda}.$
Furthermore, as in Theorem 2, such a rule exists. Now define $\eta(\lambda)$
so that $0<\eta(\lambda)<\Pi_{\Psi}(B_{\lambda}(\psi_{\text{LRSE}}(x))).$ Note
that $\eta(\lambda)\rightarrow0$ as $\lambda\rightarrow0.$ We have that, as
$\lambda\rightarrow0,$
\begin{align}
\frac{\pi_{\Psi,\lambda}(\psi_{\lambda}(\psi_{\text{LRSE}}(x))\,|\,x)}%
{\max(\eta(\lambda),\pi_{\Psi,\lambda}(\psi_{\lambda}(\psi_{\text{LRSE}%
}(x)))}  &  =\frac{\pi_{\Psi,\lambda}(\psi_{\lambda}(\psi_{\text{LRSE}%
}(x))\,|\,x)}{\pi_{\Psi,\lambda}(\psi_{\lambda}(\psi_{\text{LRSE}}%
(x)))}\nonumber\\
&  \rightarrow\frac{\pi_{\Psi}(\psi_{\text{LRSE}}(x)\,|\,x)}{\pi_{\Psi}%
(\psi_{\text{LRSE}}(x))}. \label{eq7a}%
\end{align}

Let $\epsilon>0.$ Let $\lambda_{0}$ be such that $\sup_{\psi\in\Psi}%
$diam$(B_{\lambda}(\psi))<\epsilon/2$ for all $\lambda<\lambda_{0}.$ Then for
$\lambda<\lambda_{0},$ and any $\delta(x)$ satisfying $||\delta(x)-\psi
_{\text{LRSE}}(x)||\geq\epsilon,$ we have
\begin{align}
&  \frac{\pi_{\Psi,\lambda}(\psi_{\lambda}(\delta(x))\,|\,x)}{\pi
_{\Psi,\lambda}(\psi_{\lambda}(\delta(x)))}=\frac{\int_{B_{\lambda}%
(\psi_{\lambda}(\delta(x)))}\pi_{\Psi}(\psi\,|\,x)\,\nu_{\Psi}(d\psi)}%
{\int_{B_{\lambda}(\psi_{\lambda}(\delta(x)))}\pi_{\Psi}(\psi)\,\nu_{\Psi
}(d\psi)}\nonumber\\
&  =\frac{\int_{B_{\lambda}(\psi_{\lambda}(\delta(x)))}\frac{\pi_{\Psi}%
(\psi\,|\,x)}{\pi_{\Psi}(\psi)}\pi_{\Psi}(\psi)\,\nu_{\Psi}(d\psi)}%
{\int_{B_{\lambda}(\psi_{\lambda}(\delta(x)))}\pi_{\Psi}(\psi)\,\nu_{\Psi
}(d\psi)}\nonumber\\
&  \leq\sup_{\{\psi:||\psi-\psi_{\text{LRSE}}(x)||>\epsilon/2\}}\frac
{\pi_{\Psi}(\psi\,|\,x)}{\pi_{\Psi}(\psi)}<\frac{\pi_{\Psi}(\psi_{\text{LRSE}%
}(x)\,|\,x)}{\pi_{\Psi}(\psi_{\text{LRSE}}(x))}. \label{eq7b}%
\end{align}
\ By (\ref{eq7a}) and (\ref{eq7b}) there exists $\lambda_{1}<\lambda_{0}$ such
that, for all\ $\lambda<\lambda_{1},$%
\begin{equation}
\frac{\pi_{\Psi,\lambda}(\psi_{\lambda}(\psi_{\text{LRSE}}(x))\,|\,x)}%
{\pi_{\Psi,\lambda}(\psi_{\lambda}(\psi_{\text{LRSE}}(x)))}>\sup
_{\{\psi:||\psi-\psi_{\text{LRSE}}(x)||>\epsilon/2\}}\frac{\pi_{\Psi}%
(\psi\,|\,x)}{\pi_{\Psi}(\psi)}. \label{eq7c}%
\end{equation}
Therefore, when $\lambda<\lambda_{1},$ a Bayes rule $\delta_{\lambda
,\eta(\lambda)}(x)$ satisfies
\begin{align}
&  \frac{\pi_{\Psi,\lambda}(\delta_{\lambda,\eta(\lambda)}(x)\,|\,x)}%
{\pi_{\Psi,\lambda}(\delta_{\lambda,\eta(\lambda)}(x))}\geq\frac{\pi
_{\Psi,\lambda}(\delta_{\lambda,\eta(\lambda)}(x)\,|\,x)}{\max(\eta
(\lambda),\pi_{\Psi,\lambda}(\delta_{\lambda,\eta(\lambda)}(x)))}\nonumber\\
&  \geq\frac{\pi_{\Psi,\lambda}(\psi_{\lambda}(\psi_{\text{LRSE}}%
(x))\,|\,x)}{\max(\eta(\lambda),\pi_{\Psi,\lambda}(\psi_{\lambda}%
(\psi_{\text{LRSE}}(x)))}=\frac{\pi_{\Psi,\lambda}(\psi_{\lambda}%
(\psi_{\text{LRSE}}(x))\,|\,x)}{\pi_{\Psi,\lambda}(\psi_{\lambda}%
(\psi_{\text{LRSE}}(x)))}. \label{eq7d}%
\end{align}
By (\ref{eq7b}), (\ref{eq7c}) and (\ref{eq7d}) this implies that
$||\delta_{\lambda,\eta(\lambda)}-\psi_{\text{LRSE}}(x)||<\epsilon$ and the
convergence is established.\medskip

\noindent\textbf{Proof of Corollary 7}: Following the proof of Theorem 6 we
have that\newline$\pi_{\Psi,\lambda}(\hat{\psi}_{\lambda}(x)\,|\,x)/\pi
_{\Psi,\lambda}(\hat{\psi}_{\lambda}(x))\geq\pi_{\Psi,\lambda}(\delta
_{\lambda,\eta(\lambda)}(x)\,|\,x)/\pi_{\Psi,\lambda}(\delta_{\lambda
,\eta(\lambda)}(x))$ and so by (\ref{eq7b}), (\ref{eq7c}) and (\ref{eq7d})
this implies that $||\hat{\psi}_{\lambda}(x)-\psi_{\text{LRSE}}(x)||<\epsilon$
and the convergence of $\hat{\psi}_{\lambda}(x)\ $to $\psi_{\text{LRSE}}(x)$
is established.\medskip

\noindent\textbf{Proof of Theorem 9}: For $c>0$ let $S_{c}(x)=\{\pi_{\Psi
}(\psi\,|\,x)/\pi_{\Psi}(\psi)\geq c\}$ and $S_{\eta,c}(x)=\{\pi_{\Psi}%
(\psi\,|\,x)/\max(\eta,\pi_{\Psi}(\psi))\geq c\}.$ Note that $S_{\eta
,c}(x)\uparrow S_{c}(x)$ as $\eta\rightarrow0.$

Suppose $c$ is such that $\Pi_{\Psi}(S_{c}(x)\,|\,x)\leq\gamma.$ Then
$\Pi_{\Psi}(S_{\eta,c}(x)\,|\,x)\leq\gamma$ for all $\eta$ and so $S_{\eta
,c}(x)\subset L_{\eta,\gamma}(x).$ This implies that $S_{c}(x)\subset\lim
\inf_{\eta\rightarrow0}L_{\eta,\gamma}(x)$ and since $\Pi_{\Psi}(C_{\gamma
}(x)\,|\,x)=\gamma$ this implies that $C_{\gamma}(x)\subset\lim\inf
_{\eta\rightarrow0}L_{\eta,\gamma}(x).$

Now suppose $c$ is such that $\Pi_{\Psi}(S_{c}(x)\,|\,x)>\gamma.$ Then there
exists $\eta_{0}$ such that for all $\eta<\eta_{0}$ we have $\Pi_{\Psi
}(S_{\eta,c}(x)\,|\,x)>\gamma.$ Since $L_{\eta,\gamma}(x)\subset S_{\eta
,c}(x)$ we have that $\lim\sup_{\eta\rightarrow0}L_{\eta,\gamma}(x)\subset
S_{c}(x).$ Then choosing $c=c_{\gamma^{\prime}}(x)$ for $\gamma^{\prime
}>\gamma$ implies that $\lim\sup_{\eta\rightarrow0}L_{\eta,\gamma}(x)\subset
C_{\gamma^{\prime}}(x).$\medskip

\noindent\textbf{Proof of Theorem 10}: Let $S_{c}(x)=\{\psi:\pi_{\Psi}%
(\psi\,|\,x)/\pi_{\Psi}(\psi)\geq c\}$ and $S_{\lambda,c}(x)=\{\psi:\Pi_{\Psi
}(B_{\lambda}(\psi)\,|\,x)/\Pi_{\Psi}(B_{\lambda}(\psi))\geq c\}.$ Recall
that
\[
\lim_{\lambda\rightarrow0}\Pi_{\Psi}(B_{\lambda}(\psi)\,|\,x)/\Pi_{\Psi
}(B_{\lambda}(\psi))=\pi_{\Psi}(\psi\,|\,x)/\pi_{\Psi}(\psi)
\]
for every $\psi.$ If $\pi_{\Psi}(\psi\,|\,x)/\pi_{\Psi}(\psi)>c,$ we have that
there exists $\lambda_{0}$ such that for all $\lambda<\lambda_{0},$ then
$\Pi_{\Psi}(B_{\lambda}(\psi)\,|\,x)/\Pi_{\Psi}(B_{\lambda}(\psi))>c$ and this
implies that $\psi\in\lim\inf_{\lambda\rightarrow0}S_{\lambda,c}(x).$ Now
$\Pi_{\Psi}(\pi_{\Psi}(\psi\,|\,x)/\pi_{\Psi}(\psi)=c)=0$ and so we have
$S_{c}(x)\subset\lim\inf_{\lambda\rightarrow0}S_{\lambda,c}(x)$ (after
possibly deleting a set of $\Pi_{\Psi}$-measure 0 from $S_{c}(x)).$ Now, if
$\psi\in\lim\sup_{\lambda\rightarrow0}S_{\lambda,c}(x),$ then $\Pi_{\Psi
}(B_{\lambda}(\psi)\,|\,x)/\Pi_{\Psi}(B_{\lambda}(\psi))$\newline$\geq c$ for
infinitely many $\lambda\rightarrow0,$ which implies that $\pi_{\Psi}%
(\psi\,|\,x)/\pi_{\Psi}(\psi)\geq c,$ and therefore $\psi\in$ $S_{c}(x).$ This
proves $S_{c}(x)=\lim_{\lambda\rightarrow0}S_{\lambda,c}(x)$ (up to a set of
$\Pi_{\Psi}$-measure 0) so that $\lim_{\lambda\rightarrow0}\Pi_{\Psi
}(S_{\lambda,c}(x)\Delta S_{c}(x)\,|\,x)=0$ for any $c.$

Let $c_{\lambda,\gamma}(x)=\sup\{c\geq0:\Pi_{\Psi}(S_{\lambda,c}%
(x)\,|\,x)\geq\gamma\}$ so $S_{c_{\gamma}(x)}(x)=C_{\gamma}(x)$ and
$S_{\lambda,c_{\lambda,\gamma}(x)}(x)=C_{\lambda,\gamma}(x).$ Then we have
that
\begin{align}
&  \Pi_{\Psi}(C_{\gamma}(x)\Delta C_{\lambda,\gamma}(x)\,|\,x)=\Pi_{\Psi
}(S_{c_{\gamma}(x)}(x)\Delta S_{\lambda,c_{\lambda,\gamma}(x)}%
(x)\,|\,x)\nonumber\\
&  \leq\Pi_{\Psi}(S_{c_{\gamma}(x)}(x)\Delta S_{\lambda,c_{\gamma}%
(x)}(x)\,|\,x)+\Pi_{\Psi}(S_{\lambda,c_{\lambda,\gamma}(x)}(x)\Delta
S_{\lambda,c_{\gamma}(x)}(x)\,|\,x).\label{ineq}%
\end{align}
Since $S_{c_{\gamma}(x)}(x)=\lim_{\lambda\rightarrow0}S_{\lambda,c_{\gamma
}(x)}(x)$ we have $\Pi_{\Psi}(S_{c_{\gamma}(x)}(x)\Delta S_{\lambda,c_{\gamma
}(x)}(x)\,|\,x)\rightarrow0$ and $\Pi_{\Psi}(S_{\lambda,c_{\gamma}%
(x)}(x)\,|\,x)\rightarrow\Pi_{\Psi}(S_{c_{\gamma}(x)}(x)\,|\,x)=\gamma$ as
$\lambda\rightarrow0.$ Now consider the second term in (\ref{ineq}). Since
$\pi_{\Psi}(\psi\,|\,x)/\pi_{\Psi}(\psi)$ has a continuous posterior
distribution, we have $\Pi_{\Psi}(\pi_{\Psi}(\psi\,|\,x)/\pi_{\Psi}(\psi)\geq
c\,|\,x)$ is continuous in $c.$ Let $\epsilon>0$ and note that for all
$\lambda$ small enough, $\Pi_{\Psi}(S_{\lambda,c_{\gamma-\epsilon}%
(x)}(x)\,|\,x)<\gamma$ and $\Pi_{\Psi}(S_{\lambda,c_{\gamma+\epsilon}%
(x)}(x)\,|\,x)>\gamma$ which implies that $c_{\gamma+\epsilon}(x)\leq
c_{\lambda,\gamma}(x)\leq c_{\gamma-\epsilon}(x)$ and therefore $S_{\lambda
,c_{\gamma+\epsilon}(x)}(x)\subset S_{\lambda,c_{\lambda,\gamma}(x)}\subset
S_{\lambda,c_{\gamma-\epsilon}(x)}(x).$ As $S_{\lambda,c_{\lambda,\gamma}%
(x)}(x)\subset S_{\lambda,c_{\gamma}(x)}(x)$ or $S_{\lambda,c_{\lambda,\gamma
}(x)}(x)\supset S_{\lambda,c_{\gamma}(x)}(x)$ then
\[
\Pi_{\Psi}(S_{\lambda,c_{\lambda,\gamma}(x)}(x)\Delta S_{\lambda,c_{\gamma
}(x)}(x)\,|\,x)=|\Pi_{\Psi}(S_{\lambda,c_{\lambda,\gamma}(x)}(x)\,|\,x)-\Pi
_{\Psi}(S_{\lambda,c_{\gamma}(x)}(x)\,|\,x)|.
\]
For all $\lambda$ small $|\Pi_{\Psi}(S_{\lambda,c_{\lambda,\gamma}%
(x)}(x)\,|\,x)-\Pi_{\Psi}(S_{\lambda,c_{\gamma}(x)}(x)\,|\,x)|$ is bounded
above by
\begin{align*}
& \max\{|\Pi_{\Psi}(S_{\lambda,c_{\gamma+\epsilon}(x)}(x)\,|\,x)-\Pi_{\Psi
}(S_{\lambda,c_{\gamma}(x)}(x)\,|\,x)|,\\
& \qquad|\Pi_{\Psi}(S_{\lambda,c_{\gamma-\epsilon}(x)}(x)\,|\,x)-\Pi_{\Psi
}(S_{\lambda,c_{\gamma}(x)}(x)\,|\,x)|\}
\end{align*}
and this upper bound converges to $\epsilon$ as $\lambda\rightarrow0.$ Since
$\epsilon$ is arbitrary we have that the second term in (\ref{ineq}) goes to 0
as $\lambda\rightarrow0$ and this proves the result.\medskip

\noindent\textbf{Proof of Theorem 11}: Suppose, without loss of generality
that $0<\gamma<1.$ Let $\epsilon>0$ and $\delta>0$ satisfy $\gamma+\delta
\leq1.$ Put $\gamma^{\prime}(\lambda,\gamma)=\Pi_{\Psi}(C_{\lambda,\gamma
}(x)\,|\,x),\gamma^{\prime\prime}(\lambda,\gamma)$\newline$=\Pi_{\Psi
}(C_{\gamma+\delta}(x)\,|\,x)$ and note that $\gamma^{\prime}(\lambda
,\gamma)\geq\gamma,\gamma^{\prime\prime}(\lambda,\gamma)\geq\gamma+\delta.$ By
Theorem 10 we have that $C_{\lambda,\gamma}(x)\rightarrow C_{\gamma}(x)$ and
$C_{\lambda,\gamma+\delta}(x)\rightarrow C_{\gamma+\delta}(x)$ as
$\lambda\rightarrow0$ so $\gamma^{\prime}(\lambda,\gamma)\rightarrow\gamma$
and $\gamma^{\prime\prime}(\lambda,\gamma)\rightarrow\gamma+\delta$ as
$\lambda\rightarrow0.$ This implies that there is a $\lambda_{0}(\delta)$ such
that for all $\lambda<\lambda_{0}(\delta)$ then $\gamma^{\prime}%
(\lambda,\gamma)<\gamma^{\prime\prime}(\lambda,\gamma).$ Therefore, by Theorem
9, we have that for all $\lambda<\lambda_{0}(\delta)$%
\begin{equation}
C_{\lambda,\gamma}(x)\subset\underset{\eta\rightarrow0}{\lim\inf\,}%
L_{\eta,\lambda,\gamma^{\prime}(\lambda,\gamma)}(x)\subset\underset{\eta
\rightarrow0}{\lim\sup\,}L_{\eta,\lambda,\gamma^{\prime}(\lambda,\gamma
)}(x)\subset C_{\lambda,\gamma+\delta}(x).\label{set1}%
\end{equation}
From (\ref{set1}) and Theorem 10 we have that $C_{\gamma}(x)\subset
\underset{\lambda\rightarrow0}{\lim\inf\,}\underset{\eta\rightarrow0}{\lim
\inf\,}L_{\eta,\lambda,\gamma^{\prime}(\lambda,\gamma)}(x)$\newline%
$\subset\underset{\lambda\rightarrow0}{\lim\sup}\underset{\eta\rightarrow
0}{\,\lim\sup\,}L_{\eta,\lambda,\gamma^{\prime}(\lambda,\gamma)}(x)\subset
C_{\gamma+\delta}(x).$ Since $\lim_{\delta\rightarrow0}C_{\gamma+\delta
}(x)=C_{\gamma}(x)$ this establishes the result.\newpage

\begin{center}
{\LARGE References}
\end{center}

\begin{description}
\item[ ] Aitkin, M. (2010). Statistical Inference: An Integrated
Bayesian/Likelihood Approach. Chapman and Hall/CRC, Boca Raron, FL.

\item[ ] Baldi, P. and Itti, L. (2010). Of bits and wows: A Bayesian theory of
surprise with applications to attention. Neural Networks, 23, 649-666.

\item[ ] Berger, J.O. (1985). Statistical Decision Theory and Bayesian
Analysis. \newline Springer, New York.

\item[ ] Berger, J.O., Liseo, B. and Wolpert, R.L. (1999). Integrated
likelihood methods for eliminating nuisance parameters. Stat. Sci., 14(1):1-28.

\item[ ] Bernardo, J.~M. (2005). Intrinsic credible regions: an objective
{B}ayesian approach to interval estimation. Test, 14(2):317--384. With
comments and a rejoinder by the author.

\item[ ] Bernardo, J.~M. and Smith, A. F.~M. (2000). Bayesian Theory. Wiley
Series in Probability and Statistics. John Wiley \& Sons Ltd., New York. Paperback.

\item[ ] Bishop, C.M. (2006). Pattern Recognition and Machine Learning.
Springer, New York.

\item[ ] Dempster, A. P. (1973). The direct use of likelihood for significance
testing. Memoirs, No. 1, Proceedings of Conference on Foundational Questions
in\ Statistical\ Inference, eds. O. Barndorff-Nielsen, P. Blaesild and G.
Schou, Institute of Mathematics, U. of Aarhus, 335-354 (reprinted in
Statistics and Computing (1997), 7, 247-252).

\item[ ] Evans, M. (1997). Bayesian inference procedures derived via the
concept of relative surprise. Comm. Statist. Theory Methods, 26(5):1125--1143.

\item[ ] Evans, M.~J., Guttman, I., and Swartz, T. (2006). Optimality and
computations for relative surprise inferences. Canad. J. Statist, 34(1):113-129.

\item[ ] Evans, M. and Shakhatreh, M. (2008). Optimal properties of some
{B}ayesian inferences. Electron. J. Stat., 2, 1268--1280.

\item[ ] Jang, G.~H. (2010). Invariant Procedures for Model Checking, Checking
for Prior-Data Conflict and Bayesian Inference. Ph.D. thesis, University of Toronto.

\item[ ] Le~Cam, L. (1953). On some asymptotic properties of maximum
likelihood estimates and related {B}ayes' estimates. Univ. California Publ.
Statist, 1, 277--329.

\item[ ] Robert, C.P. (1996). Intrinsic losses. Theory and Decision, 40, 191-214.

\item[ ] Rudin, W. (1974). Real and Complex Analysis. McGraw Hill, New York.

\item[ ] Tjur, T. (1974). Conditional Probability Distributions. Institute of
Mathematical\ Statistics, University of Copenhagen, Copenhagen.
\end{description}

\end{document}